\documentclass{amsart}
\usepackage{amssymb,amsthm,amsmath}   
\usepackage{aspmproc}  
\def\Left{}
\def\Right{}
\def\Center{}

\newtheorem{theorem}{Theorem}

\title[A metric discrepancy result for geometric progression with ratio 3/2]
{A metric discrepancy result for geometric progression with ratio 3/2}
\author[K. Fukuyama]{
Katusi Fukuyama
     }    
\address{Department of Mathematics\\ Kobe University 
\\ Rokko Kobe 657-8501\\ Japan}

\thanks{
The research is partially supported by JSPS KAKENHI 16K05204 and 15KT0106.
}
\dedicatory{Dedicated to Professor Norio K\^ono on his 80th birthday.}
\begin{document}
\maketitle
\begin{abstract}
We prove the law of the iterated logarithm for discrepancies of 
the sequence $\{(3/2)^k x\}$. 
\end{abstract}

\section{Introduction}

A sequence $\{x_k\}$ of real numbers is said to be uniformly distributed mod 1 (u.d. mod 1) if
$$
\lim_{N\to\infty}
\frac1N {}^\# \{ k\le N \mid \langle x_k \rangle \in [\,a, b)\}
= b-a\quad\hbox{for all}\quad0\le a< b\le 1,
$$
where $\langle x\rangle$ denotes the fractional part $x-[\,x\,]$ of $x$. 
Since the convergence is uniform in $a$ and $b$, 
we use the following discrepancy $D_N(\{x_k\})$ to 
measure the speed of convergence: 
$$
D_N(\{x_k\})
= \sup _{0\le a < b\le 1} 
\biggl|
\frac1N {}^\# \{ k\le N \mid \langle x_k \rangle \in [\,a, b)\}
- (b-a )\biggr|. 
$$

For a sequence $\{n_k\}$ with $n_{k+1}-n_k \ge C> 0$, Weyl \cite{1916W} proved 
$D_N(\{n_k x\})\to 0$ for almost every $x$. 
Erd\H os-G\'al conjecture (Cf. \cite{1964NI}) stated that 
$\{n_k x\}$ imitates the behavior of independent random variables and 
satisfies $ND_N\{n_k x\}=O((N\log\log N)^{1/2})$ 
under Hadamard's gap condition $n_{k+1}/ n_k \ge q > 1$. 
Philipp \cite{1975Ph} applied the method of Takahashi \cite{1962T} 
and solved the conjecture affirmatively by proving
$$
\frac1{4\sqrt 2}
\le
\varlimsup _{N\to \infty}
\frac{ND_N (\{n_k x\})}{\sqrt {2N\log \log N}}
\le 
\frac1{\sqrt 2}
\left (166 + \frac{664}{q^{1/2} -1}\right)
\quad\hbox{a.e.}
$$
Recently the upper bound  was replaced by 
$(1/2)\bigl({1+{4/\sqrt3}({q-1})}\bigr)^{1/2}$  (Cf. \cite{2012FM})
and 
 the lower bound  by $1/(2\sqrt2)$
(Cf.  \cite{2015AF}). 

For geometric progressions $\{\theta^k x\}$ with  $|\theta|>1$, 
we can prove the law of the iterated logarithm in exact form as below (\cite{2008F,geom}):
$$
\varlimsup _{N\to \infty}
\frac{ND_N (\{\theta^k x\})}{\sqrt {2N\log \log N}}
=
\Sigma_\theta
\quad \hbox{a.e. } x.
$$
Here $\Sigma_\theta\ge 1/2$ is a constant determined by $\theta$ 
in the following way. 

If $\theta ^k \notin {\mathbf Q}$ for all $k=1$, $2$, \dots, 
then 
$$
\Sigma_\theta=\frac12.
$$
Otherwise denote by $r$ the smallest positive $k$ with $\theta^k\in {\bf Q}$, and 
denote $\theta^r = p/q$ by $p\in \mathbf Z$ and $q\in \mathbf N$ with $\gcd(p,q)=1$. 
We can prove 
that the value $\Sigma_\theta$ does not depend on $r$ and is determined only by $p$ and $q$, i.e., 
$\Sigma_\theta=\Sigma_{p/q}$. 

When both of $p$ and $q$ are odd numbers, then (\cite{2008F,geom})
$$
\Sigma_{p/q} = \frac12\sqrt{\frac{|p|q+1}{|p|q-1}}.
$$

If $p$ is odd, $q$ is even and $|p|/q\ge 9/4$, 
or if $p$ is even, $q$ is odd and $|p|/q\ge 4$, then (\cite{ymst})
\begin{equation}\label{Eq:large}
\Sigma_{p/q}=
\sqrt{
\begin{aligned}
&\frac{(|p|q)^{I}+1}{(|p|q)^{I}-1}
v\Bigl(\frac{|p|-q-1}{2(|p|-q)}\Bigr)
\\&
+
\frac{2(|p|q)^{I}}{(|p|q)^{I}-1}
\sum_{m=1}^{I-1}
\frac1{(|p|q)^{m}}
v\Bigl(q^m\frac{|p|-q-1}{2(|p|-q)}\Bigr)
\end{aligned}
}\quad,
\end{equation}
where 
$I=\min\{ n\in {\bf N}\mid q^n = \pm 1 \hbox{ \rm mod } |p|-q\}$ and
$v(x)=\langle x \rangle(1-\langle x\rangle)$.

In particular, when $q=1$, we have
$$
\Sigma_p=
\frac12\sqrt{\frac{|p|+1}{|p|-1}} \quad \hbox{or}\quad
\Sigma_p=
\frac12\sqrt{ \frac{(|p|+1)|p|(|p|-2)}{(|p|-1)^3}}
$$
according as
$p$ is odd or  $|p|\ge 4$ is even. These cases are included in the above formula (\ref{Eq:large}), 
but the proof should be done separately (\cite{2008F,geom}). 

As  excluded examples of the fomula (\ref{Eq:large}), we have (\cite{2008F,mtwo})
$$
\Sigma_2= \frac19{\sqrt{42}}
\quad\hbox{and}\quad
\Sigma_{-2}= \frac1{49}{\sqrt{910}}
.$$

From now on we restrict ourselves to the case of positive $p$. 
If $p$ is positive,  we can prove (\cite{2008F}), 
\begin{align*}
\Sigma_{\theta}^2
&=
\sup_{0\le a \le 1}
\sigma_{\theta}^2(a)
=
\sup_{0\le a \le 1}
\biggl(
V(a ,a)
+
2\sum_{k=1}^\infty
\frac{
V(\langle p^k a \rangle,   \langle q^k a \rangle)}
{p^kq^k}
\biggr),
\end{align*}
where $V(x,\xi)= x\wedge \xi - x\xi$. 

When $pq$ is even and $p/q$ is large, the formula (\ref{Eq:large}) derived from
\begin{equation}\label{Eq:typei}
\Sigma_{p/q} = \sigma_{p/q}\Bigl(\frac{(p-q-1)/2}{p - q}\Bigr).
\end{equation}

There is a tendency that $V(\langle p^k a \rangle,   \langle q^k a \rangle)$ has 
large local maximum when $\langle p^k a \rangle=   \langle q^k a \rangle$ or
$a=n/(p^k - q^k)$ holds for some  $n\in \mathbf N$. 
Actually we can find several $p/q$ such that
\begin{equation}\label{Eq:typek}
\Sigma_{p/q} = \sigma_{p/q}\Bigl(\frac{n}{p^k - q^k}\Bigr) 
\end{equation}
holds for some $n\in {\mathbf N}$, 
and we say that $p/q$ is of type $k$. 
The formula (\ref{Eq:typei}) shows that $p/q$ is of type I if $pq$ is even and $p/q$ is large. 
Since we have
$$\Sigma_{2} = \sigma_{2}\Bigl(\frac{1}{2^2 - 1^2}\Bigr), $$
we see that $2=2/1$ is of type II. 

We are now in a position to state our result. 
It asserts that $3/2$ is of type VI. 
\begin{theorem}
We have
\begin{align*}
\Sigma_{3/2}&= \sigma_{3/2}\Bigl(\frac{277}{3^6-2^6}\Bigr)
=
\frac{2}{665}\sqrt{\frac{305671451762616889661445636790873}{10314424798490535546171949055}}.
\end{align*}
\end{theorem}

Before closing introduction, 
we state related results. 

Although $13/6$ is less than the threshold, 
we can still prove that it is of type I: 
$$
\Sigma_{13/6}=\sigma_{13/6}\Bigl(\frac{3}{13-6}\Bigr)=
\frac27\sqrt{\frac{237}{77}}.
$$

We can also find numbers of types II, III, IV, and V:
\begin{align*}
\Sigma_{4/3}&=\sigma_{4/3}\Bigl(\frac{3}{4^2-3^2}\Bigr),&
\quad
\Sigma_{8/3}&=\sigma_{8/3}\Bigl(\frac{24}{8^2-3^2}\Bigr),&
\\
\Sigma_{10/3}&=\sigma_{10/3}\Bigl(\frac{40}{10^2-3^2}\Bigr),&
\quad
\Sigma_{12/5}&=\sigma_{12/5}\Bigl(\frac{55}{12^2-5^2}\Bigr),&
\\
\Sigma_{17/8}&=\sigma_{17/8}\Bigl(\frac{101}{17^2-8^2}\Bigr),&
\Sigma_{19/10}&=\sigma_{19/10}\Bigl(\frac{2879}{19^3-10^3}\Bigr),&
\\
\Sigma_{12/7}&= \sigma_{12/7}\Bigl(\frac{8717}{12^4-7^4}\Bigr),&
\quad
\Sigma_{8/5}&= \sigma_{8/5}\Bigl(\frac{13690}{8^5-5^5}\Bigr).&
\end{align*}

These results will be proved in a separate paper \cite{martina}. 

We do not know if there exists number of type $k$ for all $k$. 

\section{Preliminary}
We denote $\sigma_{3/2}^2(x)$ simply by $\sigma^2(x)$
and put
$\displaystyle
c = \frac{277}{665}= \frac{277}{3^6-2^6}
$.
By $2^{36}=3^{36}=1$ mod 665, 
we see
$\langle 2^{k+36j}c\rangle
= \langle 2^{k}c\rangle$, 
$\langle 3^{k+36j}c\rangle
=\langle 3^{k}c\rangle$,
and
\begin{align*}
\sigma^2(c)&=
V(c, c)
+
2\sum_{k=1}^{36} \sum_{j=0}^\infty\frac1{6^{k+36j}}
V\bigl(\langle 2^{k+36j}c\rangle, \langle 3^{k+36j}c\rangle\bigr) 
\\&=
V(c, c)
+
\frac{2\cdot 6^{36}}{6^{36}-1}\sum_{k=1}^{36}\frac1{6^{k}} 
V\bigl(\langle 2^{k}c\rangle, \langle 3^{k}c\rangle\bigr) 
\\&=
\frac{1222685807050467558645782547163492}{4561296506512477081905890170847375}
.
\end{align*}
By $0\le V(x,\xi)\le 1/4$, we can verify
\begin{align}\label{Eq:tail}
\biggl|
2\sum_{n=N}^\infty
\frac{1}{6^n}
V\bigl( \langle 2^n x\rangle, \langle 3^n x\rangle\bigr)
\biggr|
&\le 
\frac{1}{10 \cdot 6^{N-1} }.
\end{align}

If $x> y$, then we have $V(x,y) = y(1-x)$, 
$V_x(x,y) = -y<0$, and $V_y(x,y) = 1-x> 0$. 
If $x< y$, then we have $V(x,y) = x(1-y)$, 
$V_x(x,y) = 1-y>0$, and $V_y(x,y) = -x< 0$. 
Since one of  $V_x(x,y)$ and  $V_y(x,y)$ is positive and the other is negative, 
we see
$
\bigl|
\frac{d}{dx} 
V\bigl( \langle 2^n x\rangle, \langle 3^n x\rangle\bigr)
\bigr|
=
\bigl|
2^n V_x\bigl( \langle 2^n x\rangle, \langle 3^n x\rangle\bigr)
+
3^n V_y\bigl( \langle 2^n x\rangle, \langle 3^n x\rangle\bigr)
\bigr|
\le 3^n
$ 
and 
\begin{equation}\label{Eq:dtail}
\biggl|
2\sum_{n=N}^\infty
\frac{1}{6^n}
\frac{d}{dx} 
V\bigl( \langle 2^n x\rangle, \langle 3^n x\rangle\bigr)
\biggr|
\le
\frac{1}{2^{N-2}}
\quad\hbox{a.e.}
\end{equation}

We can easily prove 
$V\bigl( \langle 2^n x\rangle, \langle 3^n x\rangle\bigr)
=
\int_0^x \frac{d}{dt} V\bigl( \langle 2^n t\rangle, \langle 3^n t\rangle\bigr) \,dt
$ and
$
\sigma^2(x) - \sigma^2(x_0)
=
\int_{x_0}^x 
\bigl(
\frac{d}{dt} V\bigl( \langle  t\rangle, \langle  t\rangle\bigr) 
+2\sum_{n=1}^\infty\frac1{6^n}
\frac{d}{dt} V\bigl( \langle 2^n t\rangle, \langle 3^n t\rangle\bigr) \bigr)\,dt
$.
Hence we have
$\frac{d}{dx}\sigma^2(x) 
=
\frac{d}{dx} V\bigl( \langle  x\rangle, \langle  x\rangle\bigr) 
+2\sum_{n=1}^\infty\frac1{6^n}
\frac{d}{dx} V\bigl( \langle 2^n x\rangle, \langle 3^n x\rangle\bigr)$
and  we can conclude $\sigma^2(x)$ is increasing in $x\in (a,b)$ if 
$\frac{d}{dx}\sigma^2(x)> 0$ for a.e.  $x\in (a, b)$. 

Put $h(x) = \sigma^2(x) - \sigma^2(c)$ and $\eta =\sigma^2(c)$. 
Since we have $\sigma^2(x) = \sigma^2(1-x)$ and $\sigma^2(1/2)=1/2$ (Cf. \cite{2008F}), 
it is sufficient to prove $h(x) < 0 $ for $x\in [\,0,1/2)\setminus \{c\}$.
We divide the interval $[\,0,1/2)$ into several pieces 
$[\,0, \frac13)$, $[\,\frac13, \frac{3}{2^3})$, 
$[\, \frac{3}{2^3}, \frac25)$, $[\, \frac25, \frac{11}{3^3})$, 
$[\,\frac{11}{3^3}, \frac{87}{211})$, $[\,\frac{87}{211}, \frac{101}{3^5})$, 
$[\,\frac{101}{3^5}, \frac{857}{2059})$,
$[\,\frac{857}{2059},\frac{202}{485})$,
$[\,\frac{202}{485},\frac{7985}{19171})$,
$[\,\frac{7985}{19171},\frac{24169}{58025})$,
$[\,\frac{24169}{58025},\frac{24596}{3^{10}})$,
$[\,\frac{24596}{3^{10}},\frac{72935}{175099})$,
\break
$[\,\frac{72935}{175099},\frac{73789}{3^{11}})$,
$[\,\frac{73789}{3^{11}}, c)$, 
$[\,c, \frac{911}{3^{7}})$, 
$[\,\frac{911}{3^7}, \frac{858}{2059})$, 
$[\,\frac{858}{2059}, \frac{278}{665})$, 
$[\,\frac{278}{665}, \frac{34}{3^4})$, \break
$[\,\frac{34}{3^4}, \frac{28}{65})$,  
$[\,\frac{28}{65}, \frac{4}{3^2})$, $[\, \frac{4}{3^2}, \frac12)$, 
and prove on each. 
We introduce the following notation. 
\begin{align*}
\tau_N^2 (x)
&=
V(x, x)
+
2\sum_{k=1}^{N} \frac1{6^k}
V\bigl(\langle 2^{k}x\rangle, \langle 3^{k}x\rangle\bigr) ,
\\
g_0(x) &= x(1-x),
\\
g_1(x) &= g_0(x) +({2}/{6})(3x-1)(1-2x),
\\
g_2(x) &= g_1(x) +({2}/{6^2})(2^2x-1)(4-3^2x),
\\
g_3(x) &= g_2(x) +({2}/{6^3})(3^3x-11)(4-2^3x),
\\
g_4(x) &= g_3(x) +({2}/{6^4})(2^4x-6)(34-3^4x),
\\
g_{5}(x) &= g_4(x) +({2}/{6^5})(3^5x-101)(14-2^5x),
\\
g_{6-}(x) &= g_{5}(x) +({2}/{6^6}) (3^6x-303)(27-2^6x),
\\
g_{6+}(x) &= g_{5}(x) +({2}/{6^6})(2^6x-26)(304-3^6x),
\\
g_{7\pm}(x) &= g_{6\pm}(x) +({2}/{6^7})(2^7 x-53)(911-3^7 x),
\\
g_{8-}(x) &= g_{7-}(x) +({2}/{6^8})(2^8 x-106)(2733 - 3^8 x),
\\
g_{9-}(x) &= g_{8-}(x) +({2}/{6^9})(2^9 x-213)(8199 - 3^9 x),
\\
g_{10-}(x) &= g_{9-}(x) +({2}/{6^{10}})(3^{10}x - 24596)(427- 2^{10}x),
\\
g_{11-}(x) &= g_{10-}(x) +({2}/{6^{11}})(3^{11}x-73789)(854 - 2^{11}x).
\end{align*}

\section{\Left $[\,0, \frac13)$ part}
In this interval, $\langle 2x\rangle=2x$,  $\langle 3x\rangle = 3x$, and 
$\langle 2x\rangle\le \langle 3x\rangle$. 
We have
\begin{align*}
h(x)
&\le 
\Bigl(g_0(x) + \frac{2}{6} 2x(1-3x)\Bigr)\Bigr|_{x=\frac{5}{18}}
+ \frac1{60}- \eta
\\&
=
-\frac{4903660393458055269333329257473743}{246310011351673762422918069225758250}
<0,
\end{align*}
since the 
quadratic function 
$g_0(x) + \frac{2}{6} 2x(1-3x)$ 
has the axis  at $x=\frac{5}{18}$.

\section{\Left $[\,\frac13, \frac{3}{2^3})$ part}
In $[\,\frac13, \frac12)$, 
we see $\langle 3x\rangle = 3x-1$, $\langle 2x\rangle=2x$, and 
$\langle 3x\rangle -\langle 2x\rangle=2x= x-1< 0$ to have
\begin{equation}\label{Eq:g1}
\tau_1^2(x) = g_1(x)
\quad\hbox{for}\quad
x\in [\,\tfrac13, \tfrac12).
\end{equation}

Hence in $[\,\frac13, \frac{3}{2^3})$, we have
\begin{align*}
h(x)
&\le 
g_1(x)\Bigr|_{x=\frac{3}{2^3}}
+ \frac1{10\cdot 6}- \eta
\\&
=
-\frac{5778590326763421748314562478852239}{875768929250395599725930912802696000}
<0,
\end{align*}
since the quadratic function has the axis  at $x=\frac4{3^2} > \frac{3}{2^3}$.

\section{\Right $[\, \frac{4}{3^2}, \frac12)$ part}
We have $\langle 2^2 x \rangle = 2^2x -1$, 
$\langle 3^2 x \rangle = 3^2x -4$, 
$\langle 2^2 x \rangle -\langle 3^2 x \rangle = -5x + 3> 0$. 
By (\ref{Eq:g1}) we have
\begin{align*}
h(x)
&\le 
\Bigl(
g_1(x)
+\frac{2}{6^2} (3^2 x-4)(2-2^2x)
\Bigr)\Bigr|_{x=\frac{41}{90}}
+ \frac1{10\cdot 6^2}- \eta
\\&
=
-\frac{15967737560279378662084583483719591}{2955720136220085149075016830709099000}
<0,
\end{align*}
since the quadratic function has its axis at $x=\frac{41}{90}$. 

\section{\Left $[\, \frac{3}{2^3},  \frac{11}{3^3})$ part}

In $[\, \frac{3}{3^2}, \frac{4}{3^2})\subset [\, \frac{1}{2^2}, \frac{2}{2^2})$, 
we have
$\langle 2^2x \rangle = 2^2x -1$, 
$\langle 3^2x \rangle = 3^2x -3$, 
 $\langle 3^2x \rangle - \langle 2^2x \rangle=5x-2< 0$ if $x< \frac25$ 
and  $\langle 3^2x \rangle \ge \langle 2^2x \rangle$ otherwise. 
Therefore we can conclude that 
\begin{equation}\label{Eq:g2}
\tau_2(x) = g_2(x)\quad\hbox{for}\quad x\in[\,\tfrac25, \tfrac4{3^2}).
\end{equation}

Note that $[\, \frac{3}{2^3},  \frac{11}{3^3}) \subset [\, \frac{3}{3^2}, \frac{4}{3^2})$ 
and $[\, \frac{3}{2^3},  \frac{11}{3^3})
= [\, \frac{3}{2^3}, \frac{4}{2^3})\cap[\, \frac{10}{3^3}, \frac{11}{3^3})$. 
Hence 
in $[\,\frac{3}{2^3}, \frac{11}{3^3})$, 
we have
$\langle 2^3x \rangle = 2^3x -3$, 
$\langle 3^3x \rangle = 3^3x -10$, 
and $\langle 2^3x \rangle -\langle 3^3x \rangle = -19x+7 <0$.

Thereby in $[\,\frac{3}{2^3}, \frac{2}{5})$, by (\ref{Eq:g1}) we have
\begin{align*}
h(x)
&\le 
\Bigl(g_1(x)
+\frac{2}{6^2} (3^2 x-3)(2-2^2x)
+\frac{2}{6^3} (2^3x-3)(11-3^3x)
\Bigr)\Bigr|_{x=\frac{2}{5}}
\\&\quad
 + \frac1{10\cdot 6^3}- \eta
=
-\frac{40623864934079027775310823389107}{72980744104199633310494242733558000}
<0,
\end{align*}
since the quadratic function has the axis  at $x=\frac{91}{216}> \frac 25$.

In $[\, \frac25, \frac{11}{3^3})$, by (\ref{Eq:g2}) 
we have
\begin{align*}
h(x)
&\le 
\Bigl(g_2(x)
+\frac{2}{6^3} (2^3x-3)(11-3^3x)
\Bigr)\Bigr|_{x=\frac{607}{1512}}
 + \frac1{10\cdot 6^3}- \eta
\\&
=
-\frac{115305379420414389410460002655935827}{212811849807846130733401211811055128000}
<0,
\end{align*}
since the quadratic function has the axis  at $x=\frac{607}{1512}$.

\section{\Left $[\,\frac{11}{3^3}, \frac{87}{211})$,  \Right $[\,\frac{28}{65}, \frac{4}{3^2})$ part}

On $[\,\frac{11}{3^3}, \frac{4}{3^2})\subset [\,\frac{3}{2^3}, \frac{4}{2^3})$, 
we have
$\langle 3^3 x\rangle = 3^3x -11$, 
$\langle 2^3 x\rangle = 2^3x -3$, 
to have 
$\langle 3^3 x\rangle -\langle 2^3 x\rangle 
= 19x-8< 0$ if $x< \frac{8}{19}$ 
and $\langle 3^3 x\rangle \ge\langle 2^3 x\rangle $ otherwise. 
Hence 
\begin{equation}\label{Eq:g3}
\tau_3(x) = g_3(x) \quad\hbox{for}\quad x\in[\,\tfrac{11}{3^3}, \tfrac{8}{19}).
\end{equation}

In $[\,\frac{11}{3^3}, \frac{87}{211})\subset[\,\frac{11}{3^3}, \frac{8}{19})$,  by (\ref{Eq:g3})
we have
\begin{align*}
h(x)
&\le 
g_3(x)
\Bigr|_{x=\frac{87}{211}}
 + \frac1{10\cdot 6^3}- \eta
\\&
=
-\frac{4087562128846726808409665149171701983}{29242581374367646871548627626666621462000}
<0,
\end{align*}
since the quadratic function has the axis  at $x= \frac{317}{756} >  \frac{87}{211}$.

In $[\,\frac{28}{65}, \frac{4}{3^2})\subset[\,\frac{8}{19}, \frac{4}{3^2})$, 
by (\ref{Eq:g2}) 
we have
\begin{align*}
h(x)
&\le 
\Bigr(g_2(x)
+\frac{2}{6^3} (2^3x-3)(12-3^3x)
\Bigr)\Bigr|_{x=\frac{28}{65}}
 + \frac1{10\cdot 6^3}- \eta
\\&
=
-\frac{62145572531001959682866090587826537}{25616241180574071291983479199478858000}
<0,
\end{align*}
since the quadratic function has the axis  at $x= \frac{205}{504} < \frac{28}{65}$.

\section{\Right $[\,\frac{34}{3^4}, \frac{28}{65})$ part}
It is included in  $[\,\frac{11}{3^3}, \frac{4}{3^2})$ and 
$[\,\frac{6}{2^4}, \frac{7}{2^4})\cap [\,\frac{34}{3^4},\frac{35}{3^4})$. 
We have 
$\langle 3^4 x\rangle = 3^4 x-34$, 
$\langle 2^4 x\rangle = 2^4 x-6$, 
$\langle 3^4 x\rangle -\langle 2^4 x\rangle 
=
65x-28< 0$. 

In $[\,\frac{8}{19}, \frac{28}{65})$, by (\ref{Eq:g2}) 
we have
\begin{align*}
h(x)
&\le 
\Bigl(g_2(x)
+\frac{2}{6^3} (2^3x-3)(12-3^3x)
+\frac{2}{6^4} (3^4x-34)(7-2^4x)
\Bigr)\Bigr|_{x=\frac{8}{19}}
\\&
\qquad
 + \frac1{10\cdot 6^4}- \eta
\\&
=
-\frac{638551936297453983963801188142263}{3940960181626780198766689107612132000}
<0,
\end{align*}
since the quadratic function has the axis  at
$x=\frac{4801}{11664}< \frac{8}{19}$.

In $[\,\frac{34}{3^4},\frac{8}{19})$, by (\ref{Eq:g3}) 
we have
\begin{align*}
h(x)
&\le 
\Bigl(g_3(x)
+\frac{2}{6^4} (3^4x-34)(7-2^4x)
\Bigr)\Bigr|_{x=\frac{8}{19}}
 + \frac1{10\cdot 6^4}- \eta
\\&
=
-\frac{638551936297453983963801188142263}{3940960181626780198766689107612132000}
<0,
\end{align*}
since the quadratic function has the axis  at
$x=\frac{4915}{11664}> \frac{8}{19}$.

\section{\Left $[\,\frac{87}{211}, \frac{101}{3^5})$ part}

In $[\,\frac{87}{211}, \frac{34}{3^4})$, 
we have 
$\langle 3^4 x\rangle = 3^4 x-33$, 
$\langle 2^4 x\rangle = 2^4 x-6$, 
$\langle 3^4 x\rangle - \langle 2^4 x\rangle
=
65x-27<0$ if $x< \frac{27}{65}$ and 
$\langle 3^4 x\rangle \ge \langle 2^4 x\rangle$ otherwise. 
Hence 
\begin{equation}\label{Eq:g4}
\tau_4(x) = g_4(x) \quad\hbox{for}\quad x\in [\,\tfrac{27}{65},\tfrac{34}{3^4}).
\end{equation}

In $[\,\frac{87}{211}, \frac{101}{3^5})$, 
we see
$\langle 3^5 x\rangle = 3^5 x-100$, 
$\langle 2^5 x\rangle = 2^5 x-13$, 
$\langle 3^5 x\rangle - \langle 2^5 x\rangle
=
211x-87\ge 0$. 

Hence in $[\,\frac{87}{211}, \frac{302}{3^6})\subset [\,\frac{87}{211}, \frac{27}{65})$, 
by (\ref{Eq:g3}) 
we have
\begin{align*}
&h(x)
\\
&\le 
\Bigl(g_3(x)
+\frac{2}{6^4} (3^4x-33)(7-2^4x)
+\frac{2}{6^5} (2^5x-13)(101-3^5x)
\Bigr)\Bigr|_{x=\frac{302}{3^6}}
\\&
 + \frac1{10\cdot 6^5}- \eta
=
-\frac{13250661085528974219050820575468766883}{155139838509919829304649483410259188312000}
<0,
\end{align*}
since the quadratic function has the axis  at $x=\frac{35785}{85536}> \frac{27}{65}$.

In $[\,\frac{302}{3^6}, \frac{101}{3^5})$, 
we see
$\langle 3^6x \rangle = 3^6 x - 302$,
$\langle 2^6x \rangle = 2^6 x - 26$, 
$\langle 3^6x \rangle -\langle 2^6x \rangle =
665 x - 276 < 0$ if $x < \frac{276}{665}$
and $\langle 3^6x \rangle \ge\langle 2^6x \rangle$ otherwise. 
Note that
$\frac{302}{3^6}<  \frac{276}{665} < \frac{27}{65} < \frac{101}{3^5}$. 

In $[\,\frac{302}{3^6},  \frac{276}{665})$, by (\ref{Eq:g3}) we have
\begin{align*}
h(x)
&\le 
\Bigl(g_3(x)
+\frac{2}{6^4} (3^4x-33)(7-2^4x)
+\frac{2}{6^5} (2^5x-13)(101-3^5x)
\\&\qquad
+\frac{2}{6^6} (3^6x-302)(27-2^6 x)
\Bigr)\Bigr|_{x=\frac{276}{665}}
 + \frac1{10\cdot 6^6}- \eta
\\&
=
-\frac{9623392693609973991877849224521309}{425623699615692261466802423622110256000}
<0,
\end{align*}
since the quadratic function has the axis  
at $x=\frac{19517}{46656} > \frac{276}{665}$. 

In $[\,\frac{276}{665}, \frac{27}{65})$, by (\ref{Eq:g3}) we have
\begin{align*}
h(x)
&\le 
\Bigl(g_3(x)
+\frac{2}{6^4} (3^4x-33)(7-2^4x)
+\frac{2}{6^5} (2^5x-13)(101-3^5x)
\\&\qquad
+\frac{2}{6^6} (2^6x-26)(303-3^6x)
\Bigr)\Bigr|_{x=\frac{27}{65}}
 + \frac1{10\cdot 6^6}- \eta
\\&
=
-\frac{1906349521544493873261549971920349}{425623699615692261466802423622110256000}
<0,
\end{align*}
since the quadratic function has the axis  
at $x=\frac{1318}{3159} > \frac{27}{65}$.

In $[\,\frac{27}{65},\frac{101}{3^5})$, by (\ref{Eq:g4}) we have
\begin{align*}
h(x)
&\le 
\Bigl(
g_4(x)
+\frac{2}{6^5} (2^5x-13)(101-3^5x)
\\&\qquad
+\frac{2}{6^6} (2^6x-26)(303-3^6x)
\Bigr)\Bigr|_{x=\frac{27}{65}}
 + \frac1{10\cdot 6^6}- \eta
\\&
=
-\frac{1906349521544493873261549971920349}{425623699615692261466802423622110256000}
<0,
\end{align*}
since the quadratic function has the axis  
at $x=\frac{20893}{50544} < \frac{27}{65}$.

\section{\Right $[\,\frac{278}{665}, \frac{34}{3^4})$ part}

In $[\,\frac{101}{3^5}, \frac{34}{3^4})$, 
we have
$\langle 3^5 x\rangle = 3^5 x-101$, 
$\langle 2^5 x\rangle = 2^5 x-13$, 
$\langle 3^5 x\rangle - \langle 2^5 x\rangle =211x -88< 0$ if $x< \frac{88}{211}$
and $\langle 3^5 x\rangle\ge \langle 2^5 x\rangle $ otherwise. 
Therefore 
\begin{equation}\label{Eq:g5}
\tau_5(x) = g_5(x)\quad\hbox{for}\quad x\in[\,\tfrac{101}{3^5}, \tfrac{88}{211}).
\end{equation} 
Note that $c< \frac{88}{211}< \frac{278}{665}$. 
In $[\,  \frac{278}{665},  \frac{34}{3^4})$, we have
\begin{align*}
h(x)
&\le 
\Bigl(g_4(x)
+\frac{2}{6^5} (2^5x-13)(102-3^5x)
\Bigr)\Bigr|_{x=\frac{278}{665}}
 + \frac1{10\cdot 6^5}- \eta
\\&
=
-\frac{5091905468453476674801592843459949}{70937283269282043577800403937018376000}
<0,
\end{align*}
since the quadratic function has the axis  
at $x=\frac{11809}{28512} < \frac{278}{665}$.

\section{\Right $[\,\frac{858}{2059}, \frac{278}{665})$ part}

Note that $c< \frac{858}{2059}<  \frac{304}{3^6}
< \frac{88}{211}< \frac{278}{665}$. 

In $[\,\frac{304}{3^6}, \frac{278}{665})$, 
we see
$\langle 3^6 x\rangle = 3^6 x - 304$, 
$\langle 2^6 x\rangle = 2^6 x - 26$, 
$\langle 3^6 x\rangle - \langle 2^6 x\rangle = 665x-278< 0$.

In $[\,\frac{88}{211}, \frac{278}{665})$, by (\ref{Eq:g4}) we have
\begin{align*}
h(x)
&\le 
\Bigl(g_4(x)
+\frac{2}{6^5} (2^5x-13)(102-3^5x)
\\&\qquad\qquad
+\frac{2}{6^6} (3^6 x - 304)(27-2^6x)
\Bigr)\Bigr|_{x=\frac{88}{211}}
 + \frac1{10\cdot 6^6}- \eta
\\&
=
-\frac{56411054147338655852222279153714125703}{6316397576863411724254503567359990235792000}
<0,
\end{align*}
since the quadratic function has the axis  
at $x=\frac{251701}{606528} < \frac{88}{211}$. 

In $[\,\frac{304}{3^6},\frac{88}{211})$, by (\ref{Eq:g5}) we have
\begin{align*}
h(x)
&\le 
\Bigl(g_5(x)
+\frac{2}{6^6} (3^6 x - 304)(27-2^6x)
\Bigr)\Bigr|_{x=\frac{88}{211}}
 + \frac1{10\cdot 6^6}- \eta
\\&
=
-\frac{56411054147338655852222279153714125703}{6316397576863411724254503567359990235792000}
<0,
\end{align*}
since the quadratic function has the axis  
at $x=\frac{19459}{46656} > \frac{88}{211}$. 
                                      
In $[\,\frac{101}{3^5}, \frac{304}{3^6})$, 
we see $\langle 3^6 x\rangle = 3^6 x -303$, 
$\langle 2^6 x\rangle = 2^6 x -26$, 
$\langle 3^6 x\rangle - \langle 2^6 x\rangle 
= 665 x - 277\ge 0 $ if $x\ge c$ and 
$\langle 3^6 x\rangle < \langle 2^6 x\rangle $ if $x<  c$. 
Hence 
\begin{align}\label{Eq:g6-}
\tau_6(x) &= g_{6-}(x) \quad\hbox{for}\quad x\in[\,\tfrac{101}{3^5}, c),
\\ \label{Eq:g6+}
\tau_6(x) &= g_{6+}(x) \quad\hbox{for}\quad x\in[\, c,\tfrac{304}{3^6})
.
\end{align} 

In $[\,\frac{858}{2059}, \frac{304}{3^6})\subset [\, c,\tfrac{304}{3^6})$,  
by (\ref{Eq:g6+}) we have
\begin{align*}
h(x)
&\le 
g_{6+}(x)
\Bigr|_{x=\frac{858}{2059}}
 + \frac1{10\cdot 6^6}- \eta
\\&
=
-\frac{245734830243268354941082685853561543661}{200491509741159404926171222855543067801904000}
<0,
\end{align*}
since the quadratic function has the axis  
at $x=\frac{126119}{303264} < \frac{858}{2059}$.

\section{\Right $[\,\frac{911}{3^7}, \frac{858}{2059})$ part}
In $[\,\frac{911}{3^7}, \frac{858}{2059})$, 
we have
$\langle 3^7x\rangle= 3^7x-911$, 
$\langle 2^7x\rangle= 2^7 x-53$, 
$\langle 3^7x\rangle-\langle 2^7x\rangle = 2059x-858<0$. 

In $[\,\frac{2627}{6305}, \frac{858}{2059})\subset [\,\frac{911}{3^7}, \frac{858}{2059})$, 
by (\ref{Eq:g6+}) 
we have 
\begin{align*}
h(x)
&\le 
\Bigl(g_{6+}(x)+\frac2{6^7}(3^7x-911)(54-2^7x)\Bigr)
\Bigr|_{x=\frac{2627}{6305}}
 + \frac1{10\cdot 6^7}- \eta
\\&
=
-\frac{296152864180283626343198670531396721}{357807656810258627806425237458320688544000}
<0,
\end{align*}
since the quadratic function has the axis  
at $x=\frac{874067}{2099520} < \frac{2627}{6305}$.

In $[\,\frac{911}{3^7}, \frac{2627}{6305})$, 
we have
$\langle 3^8 x\rangle= 3^8x-2733$, 
$\langle 2^8 x\rangle= 2^8x-106$, 
$\langle 3^8 x\rangle-\langle 2^8 x\rangle = 6305x-2627< 0$. 
Hence by (\ref{Eq:g6+}) we have
\begin{align*}
h(x)
&\le 
\Bigl(g_{6+}(x)
+\frac2{6^7}(3^7x-911)(54-2^7x)
\\&\qquad\qquad\qquad
+\frac{2}{6^8} (3^8 x - 2733)(107-2^8x)
\Bigr)\Bigr|_{x=\frac{911}{3^7}}
 + \frac1{10\cdot 6^8}- \eta
\\&
=
-\frac{2903244137571860267233333557802343701}{11170068372714227709934762805538661558464000}
<0,
\end{align*}
since the quadratic function has the axis  
at $x=\frac{3963493}{9517824} < \frac{911}{3^7}$. 

\section{\Left $[\,\frac{101}{3^5}, \frac{857}{2059})$ part}
In $[\,\frac{101}{3^5}, \frac{857}{2059})\subset [\,\frac{101}{3^5}, c)$, 
by (\ref{Eq:g6-}) we have
\begin{align*}
h(x)
&\le 
g_{6-}(x)
\Bigr|_{x=\frac{857}{2059}}
 + \frac1{10\cdot 6^6}- \eta
\\&
=
-\frac{5169596011476730532177386324915893190949}{1804423587670434644335541005699887610217136000}
<0,
\end{align*}
since the  quadratic function has the axis  at $x=\frac{84301}{202176}> \frac{857}{2059}$. 

\section{\Left $[\,\frac{857}{2059},\frac{202}{485})$ part}
In $ [\,\frac{857}{2059},\frac{911}{3^7})\subset[\,\frac{101}{3^5}, \frac{304}{3^6})$, 
we have
$\langle 3^7 x\rangle = 3^7 x-910$, 
$\langle 2^7 x\rangle = 2^7 x-53$, 
$\langle 3^7 x\rangle - \langle 2^7 x\rangle =2059x-857\ge 0$.
By (\ref{Eq:g6-}) and (\ref{Eq:g6+}) we have 
\begin{align}\label{Eq:g7-}
\tau_7(x)  &= g_{7-}(x)\quad\hbox{for}\quad x\in[\,\tfrac{857}{2059},c),
\\ \label{Eq:g7+}
\tau_7(x)  &= g_{7+}(x)\quad\hbox{for}\quad x\in[\,c, \tfrac{911}{3^7}).
\end{align} 

In $[\,\frac{857}{2059},\frac{202}{485})\subset[\,\frac{857}{2059},c)$,
by (\ref{Eq:g7-}) we have
\begin{align*}
h(x)
&\le 
g_{7-}(x)
\Bigr|_{x=\frac{202}{485}}
 + \frac1{10\cdot 6^7}- \eta
\\&
=
-\frac{103698794723921201659167205767313}{3058185100942381434242950747507014432000}
<0,
\end{align*}
since the  quadratic function has the axis  at $x=\frac{1749937}{4199040}> \frac{202}{485}$.

\section{\Left $[\,\frac{202}{485},\frac{7985}{19171})$ part}

In $[\,\frac{202}{485},c)\subset
[\,\frac{857}{2059}, c)$, 
we have
$\langle 3^8 x\rangle = 3^8x- 2732$, 
$\langle 2^8 x\rangle = 2^8x- 106$, 
$\langle 3^8 x\rangle - \langle 2^8 x\rangle = 6305(x-\frac{202}{485})\ge 0$. 
By (\ref{Eq:g7-}) we see
\begin{equation}\label{Eq:g8-}
\tau_8(x) = g_{8-}(x)\quad\hbox{for}\quad x\in [\,\tfrac{202}{485},c). 
\end{equation}
Hence in $[\,\frac{202}{485},\frac{7985}{19171})\subset [\,\frac{202}{485},c)$,
 we have
\begin{align*}
h(x)
&\le 
g_{8-}(x)
\Bigr|_{x=\frac{7985}{19171}}
 + \frac1{10\cdot 6^8}- \eta
\\&
=
-\frac{47689069442943750663023590748245592189}{15599498462223969356869812416770378735394496000}
<0,
\end{align*}
since the  quadratic function has the axis  at $x=\frac{247807}{594864}> \frac{7985}{19171}$.

\section{\Left $[\,\frac{7985}{19171},\frac{24169}{58025})$ part}
In
$[\,\frac{7985}{19171}, c)
\subset
[\,\frac{202}{485},c)$, 
we have
$\langle 3^9 x\rangle = 3^9 x -8198$, 
$\langle 2^9 x\rangle = 2^9 x -213$, 
$\langle 3^9 x\rangle - \langle 2^9 x\rangle  =
19171x-7985\ge 0$. 
By (\ref{Eq:g8-}) we see
\begin{equation}\label{Eq:g9-}
\tau_9(x) = g_{9-}(x)\quad\hbox{for}\quad x\in [\,\tfrac{7985}{19171}, c).
\end{equation}
Hence in $[\,\frac{7985}{19171},\frac{24169}{58025})\subset [\,\frac{7985}{19171},c)$, we have
\begin{align*}
h(x)
&\le 
g_{9-}(x)
\Bigr|_{x=\frac{24169}{58025}}
 + \frac1{10\cdot 6^9}- \eta
\\&
=
-\frac{5902470452852577000431786200645894683457}{2476280506033531932376735578547810572039895680000}
<0,
\end{align*}
since the  quadratic function has the axis  at $x=\frac{2954029}{7091712}> \frac{24169}{58025}$.

\section{\Left $[\,\frac{24169}{58025},\frac{24596}{3^{10}})$ part}
In 
$[\,\frac{24169}{58025},\frac{24596}{3^{10}})$,
we have
$\langle 3^{10}x\rangle = 3^{10}x - 24595$, 
$\langle 2^{10}x\rangle = 2^{10}x - 426$, 
$\langle 3^{10}x\rangle -  \langle 2^{10}x\rangle = 58025x-24169\ge 0$. 

In $[\,\frac{24169}{58025},\frac{24596}{3^{10}})
\subset[\,\frac{7985}{19171},c)$, by (\ref{Eq:g9-}) we have
\begin{align*}
&h(x)
\\
&\le 
\Bigl(g_{9-}(x)
+\frac{2}{6^{10}}(2^{10}x - 426)(24596 -3^{10}x)
\Bigr)\Bigr|_{x=\frac{24169}{58025}}
 + \frac1{10\cdot 6^{10}}- \eta
\\&
=
-\frac{11899418967263802834461455793906308311731}{4952561012067063864753471157095621144079791360000}
<0,
\end{align*}
since the  quadratic function has the axis  at $x=\frac{18889067}{45349632}< 
\frac{24169}{58025}$.

\section{\Left $[\,\frac{24596}{3^{10}},\frac{72935}{175099})$ part}
In $
[\,\frac{24596}{3^{10}},c)\subset[\,\frac{7985}{19171},c)$, 
by noting $\frac{24170}{58025}> c$, 
we have
$\langle 3^{10}x\rangle = 3^{10}x - 24596$, 
$\langle 2^{10}x\rangle = 2^{10}x - 426$, 
$\langle 3^{10}x\rangle -  \langle 2^{10}x \rangle = 58025x-24170<0$. 
By (\ref{Eq:g9-}) we have
\begin{equation}\label{Eq:g10-}
\tau_{10}(x) = g_{10-}(x)\quad\hbox{for}\quad x\in[\,\tfrac{24596}{3^{10}},c).
\end{equation}

In $[\,\frac{24596}{3^{10}},\frac{72935}{175099})\subset[\,\frac{24596}{3^{10}},c)$,
 we have
\begin{align*}
&h(x)
\le 
g_{10-}(x)
\Bigr|_{x=\frac{72935}{175099}}
 + \frac1{10\cdot 6^{10}}- \eta
\\&
=
-\frac{92702352365411529614480198024779430800965409}{16912123272164340374273932618342897873352232696576000}
<0,
\end{align*}
since the  quadratic function has the axis  at $x=\frac{528952925}{1269789696}> \frac{72935}{175099}$.

\section{\Left $[\,\frac{72935}{175099},\frac{73789}{3^{11}})$ part}

In
$
[\,\frac{72935}{175099},\frac{73789}{3^{11}})\subset[\,\frac{24596}{3^{10}},c)$, 
we have
$\langle 3^{11}x\rangle =3^{11}x - 73788$, 
$\langle 2^{11}x\rangle =2^{11}x - 853$, 
$\langle 3^{11}x\rangle - \langle 2^{11}x\rangle 
=175099x-72935\ge 0 $. 
By (\ref{Eq:g10-}) we have
\begin{align*}
h(x)
&\le 
\Bigl(g_{10-}(x)
+\frac{2}{6^{11}}(2^{11}x - 853)(73789-3^{11}x)
\Bigr)\Bigr|_{x=\frac{73789}{3^{11}}}
\\&\qquad
 + \frac1{10\cdot 6^{11}}- \eta
\\&
=
-\frac{29347546710908416454761350780039732607}{65143838749669376004339536681901474208962048000}
<0,
\end{align*}
since the  quadratic function has the axis  at $x=\frac{3475943813}{8344332288}
> \frac{73789}{3^{11}}$.

\section{\Center $[\,\frac{73789}{3^{11}}, c)$, 
$[\,c, \frac{911}{3^{7}})$
part}

In $[\,c, \frac{911}{3^{7}})$ we have
\begin{align*}
\frac{d}{dx}
\sigma^2(x)
& \le
\frac{d}{dx}
g_{7+}(x)
+ \frac1{2^6}
=
\frac{1745947}{139968}-30x
+ \frac1{2^6}
\\&
\le
\frac{1745947}{139968}-30c
+ \frac1{2^6}
=
-\frac{62497}{9307872}< 0
\quad\hbox{a.e.}
\end{align*}

In $[\,\frac{73789}{3^{11}}, c)$
we have
$\langle 3^{11}x\rangle = 3^{11}x-73789$,
$\langle 2^{11}x\rangle =  2^{11}x-853$,
$\langle 3^{11}x\rangle - \langle 2^{11}x\rangle = 175099 x - 72936< 0$, 
and 
$\tau_{11}(x) = g_{11-}(x)$.
Hence 
\begin{align*}
\frac{d}{dx}
\sigma^2(x)
& \ge
\frac{d}{dx}
g_{11-}(x)
- \frac1{2^{10}}
=
\frac{27157195}{1417176}-46x
- \frac1{2^{10}}
\\&
\ge
\frac{27157195}{1417176}-46c
- \frac1{2^{10}}
=
\frac{122591869}{120630021120}> 0
\quad\hbox{a.e.}
\end{align*}


\end{document}